\documentclass[12pt]{article}
\usepackage{amsmath,amsfonts,amssymb,amsthm,enumerate}

\newcommand{\R}{{\mathbb R}}
\newcommand{\C}{{\mathbb C}}
\newcommand{\Z}{{\mathbb Z}}
\newcommand{\Q}{{\mathbb Q}}

\newcommand{\F}{{\mathbb F}}

\newcommand{\BW}{\Lambda_{16}}
\newcommand{\ch}{{\rm ch }}

\newcommand{\Aut}{{\rm Aut}}
\newcommand{\AutO}{O}

\newcommand{\Hom}{{\rm Hom}}

\newcommand{\Cent}{Z}

\newcommand{\Ker}{{\rm Ker\ }}
\newtheorem{dfn}{Definition}[section]
\newtheorem{pro}[dfn]{Proposition}
\newtheorem{thm}[dfn]{Theorem}

\newtheorem{lem}[dfn]{Lemma}
\newtheorem{cor}[dfn]{Corollary}
\newtheorem{rem}[dfn]{Remark}
\newtheorem{note}[dfn]{Note}

\newcommand{\qe}{\qed\vskip2ex}

\makeatletter
\@addtoreset{equation}{section}

\makeatother

\def\bl{\begin{lem}}
\def\el{\end{lem}}
\def\bt{\begin{thm}}
\def\et{\end{thm}}
\def\bp{\begin{pro}}
\def\ep{\end{pro}}
\def\br{\begin{rem}\rm }
\def\er{\end{rem}}
\def\bc{\begin{cor}}
\def\ec{\end{cor}}
\def\bd{\begin{dfn}}
\def\ed{\end{dfn}}
\def\bn{\begin{note}\rm }
\def\en{\end{note}}
\def\proof{{\it Proof.}}

\title{\begin{flushright}
\end{flushright}\Large Lifts of automorphisms of vertex operator algebras\\ in simple current extensions}
\author{Hiroki SHIMAKURA\footnote{The author was supported by the COE grant of Hokkaido University and JSPS Grants-in-Aid for Scientific Research No. 18740001.}
}
\date{\small\it Department of Mathematics, Hokkaido University\\
Kita 10, Nishi 8, Kita-Ku, Sapporo, Hokkaido, 060-0810, Japan.\\
{\rm e-mail: shimakura@math.sci.hokudai.ac.jp}
}
\begin{document}
\maketitle

\begin{abstract}
In this article, we study isomorphisms between simple current extensions of a simple VOA.
For example, we classify the isomorphism classes of simple current extensions of the VOAs $V_{\sqrt2E_8}^+$ and $V_{\BW}^+$, where $\BW$ is the Barnes-Wall lattice of rank $16$.
Moreover, we consider the same simple current extension and describe the normalizer of the abelian automorphism group associated with this extension.
In particular, we regard the moonshine module $V^\natural$ as simple current extensions of five subVOAs $V_L^+$ for $2$-elementary totally even lattices $L$, and describe corresponding five normalizers of elementary abelian $2$-group in the automorphism group of $V^\natural$ in terms of $V_L^+$.
By using this description, we show that three of them form a Monster amalgam.
\end{abstract}

\section*{Introduction}
In the study of vertex operator algebras (VOAs), the automorphism groups play important roles to develop relations between VOAs and other areas.
One of the most important examples of VOAs is the moonshine module constructed in \cite{FLM}, and its automorphism group is isomorphic to the Monster.
By using this realization of the Monster, Borcherds proved the famous moonshine conjecture, an interesting correspondence between the conjugacy classes of the Monster and Hauptmodule.
So it is important to determine the automorphism group of a VOA and to study groups from the view point of VOAs.
However there are relatively few known constructions of automorphisms of VOAs such as exponentials of the zero-th product of elements of $V_1$, symmetries of the fusion rules of subVOAs (cf. \cite{Mi}) and lifts of automorphisms of combinatorial objects (cf. \cite{FLM}).
Hence we would like to find more automorphisms of VOAs.

A simple current extension is a typical method of constructing new VOAs.
In some cases, this construction extends symmetries of VOAs.
For example, the moonshine module $V^\natural$ is a simple current extension of the VOA $V_{\Lambda}^+$ associated with the Leech lattice $\Lambda$, and the automorphism group of $V^\natural$ is the Monster though that of $V_{\Lambda}^+$ is induced from $\Lambda$ (\cite{Sh2}).
So, the study of the automorphism groups of simple current extensions is important to find large symmetries of VOAs.

In this article, we construct some isomorphisms between simple current extensions of a simple VOA $V(0)$ by using automorphisms of $V(0)$.
Moreover, we consider the same simple current extension $V$ of $V(0)$  graded by an abelian group $A$ and describe the normalizer of $A^*$ in the automorphism group of $V$.
So we obtain automorphisms of $V$ induced from $V(0)$.
As applications, we obtain two results.
One is the classification of the isomorphism classes of simple current extensions of $V_L^+$ when $L$ is isomorphic to $\sqrt2E_8$ or $\BW$.
Another is the description of some normalizers in the automorphism group of $V^\natural$ of elementary abelian $2$-subgroups in terms of $V_L^+$.
Furthermore, using our description, we show that three normalizers form a Monster amalgam.

\medskip

Let us explain our method.
First, we construct isomorphisms between simple current extensions, which generalizes the results in \cite{DM2} on the uniqueness of the VOA structure of a simple current extension.
Consider a simple VOA $V(0)$ and its simple current extensions $V$ and $V^\prime$.
When $V(0)$ has an automorphism which maps the isomorphism classes of irreducible $V(0)$-modules contained in $V$ to those contained in $V^\prime$, we show that it extends to an isomorphism between $V$ and $V^\prime$.
So, by using automorphisms of $V(0)$, we obtain many isomorphisms between simple current extensions of $V(0)$.
In \cite{Sh2,Sh3}, the automorphism group of the VOA $V_L^+$ is studied by using its action on the isomorphism classes of irreducible $V_L^+$-modules.
In particular, when $L$ is isomorphic to $\sqrt2E_8$ or $\BW$, $V_L^+$ has large symmetries associated with orthogonal groups.
Applying our results to $V_{\sqrt2E_8}^+$ and $V_{\BW}^+$, we can classify the isomorphism classes of these simple current extensions.

Next, we consider automorphisms of a simple current extension $V=\oplus_{\alpha\in A}V(\alpha)$ of a simple VOA $V(0)$ graded by an abelian group $A$.
Then applying the results to the same simple current extension, we obtain lifts of automorphisms of $V(0)$ preserving the set of components in the decomposition of $V$ into irreducible $V(0)$-modules.
More precisely, the normalizer $N_{\Aut(V)}(A^*)$ of $A^*$ in the automorphism group $\Aut(V)$ of $V$ is obtained by lifts of the automorphisms of $V(0)$ preserving the set of isomorphism classes of irreducible $V(0)$-modules $V(\alpha)$, $\alpha\in A$.
By similar arguments, the centralizer $C_{\Aut(V)}(A^*)$ is described in terms of the automorphism group of $V(0)$.

As an application, we consider some subVOAs of the moonshine module $V^\natural$ isomorphic to the VOA $V_L^+$ associated with a $2$-elementary totally even full sublattice $L$, which is an even lattice $L$ such that the dual lattice $L^*$ is a sublattice of $L/2$ and $\sqrt2L^*$ is even.
Then $V^\natural$ is a simple current extension of $V_L^+$ graded by an elementary abelian $2$-group.
By using the description of $\Aut(V_L^+)$ in \cite{Sh2,Sh3}, the normalizer and centralizer of the associated $2$-group in $\Aut(V^\natural)$ are described.
In particular, we consider five subVOAs of $V^\natural$ isomorphic to $V_L^+$ or $V_{L_1}^+\otimes V_{L_2}^+$ for certain $2$-elementary totally even lattices $L$, $L_1\oplus L_2$, and describe the normalizers of corresponding elementary abelian group of order $2^i$ ($i=1,2,3,5,10)$ in $\Aut(V^\natural)$.
Moreover, by using these descriptions, we check that three of them satisfy the axiom of a Monster amalgam. 
This shows that the Monster acts faithfully on $V^\natural$ without using the structure of the Griess algebra of $V^\natural$.

\medskip

The organization of this paper is as follows:
In Section 1, we recall some definitions and facts necessary in this paper.
In Section 2, we study isomorphisms between simple current extensions of VOAs.
In Section 3, we determine the isomorphism classes of simple current extensions of $V_{\sqrt2E_8}^+$ and $V_{\BW}^+$ by using the results in Section 2.
In Section 4, we describe the normalizers in the automorphism group of $V^\natural$ of some elementary abelian $2$-groups of order $2^i$ $(i=1,2,3,5,10)$ in terms of automorphism groups of $V_L^+$.
Moreover we show that the normalizers of the elementary abelian $2$-subgroups of order $2$, $2^2$ and $2^3$ form a Monster amalgam.

\medskip

Throughout this paper, we will work over the field $\C$ of complex numbers unless otherwise stated.
We denote by $\Z$ the set of integers and by $\Z_p$ the ring of integers modulo $p$.
We often identify $\Z_2$ with the field $\F_2$ of two elements.
For a group $G$ and its subgroup $H$, $N_G(H)$ and $C_G(H)$ denote the normalizer and centralizer of $H$ in $G$ respectively.
Let $\Omega_n$ denote the set $\{1,2,\dots,n\}$ for $n\in\Z_{>0}$.
We view the power set $\mathcal{P}(\Omega_n)$ of $\Omega_n$ as an $n$-dimensional vector space over $\F_2$ naturally.
For a subset $U$ of an $n$-dimensional vector space $\R^n$ over the real field $\R$ and $m\in\R$, let $U_m$ denote the set of vectors in $U$ of norm $m$.
A sublattice $U$ of a lattice $L$ is called {\it full} if the ranks of $U$ and $L$ are the same.
A subVOA $V(0)$ of $V$ is called {\it full} if the Virasoro elements of $V(0)$ and $V$ are the same.
We often identify modules of a VOA with their respective isomorphism classes.

\medskip

{\it Acknowledgments.} The author thanks Professor Atsushi Matsuo for giving helpful advice.
He also thanks Professor Satoshi Yoshiara and Professor Alexander Anatolievich Ivanov for useful comments from the viewpoint of finite group theory.

\medskip

\section{Preliminaries}\label{ChVOA}
In this section, we recall or give some definitions and facts necessary in this paper.
For details of the axiom of vertex operator algebras, see \cite{Bo,FLM}.

\subsection{Simple current extension}\label{Ss}
We start by recalling simple current extensions.
Let $V(0)$ be a simple VOA.
An irreducible $V(0)$-module $M^1$ is called a {\it simple current} if for any irreducible $V$-module $M^2$, there exists an irreducible $V(0)$-module $M^3$ such that the fusion rule $M^1\times M^2=M^3$ holds.
A simple VOA $V$ is called a {\it simple current extension} of $V(0)$ if $V$ is a direct sum of inequivalent simple current irreducible $V(0)$-modules graded by an abelian group $A$, namely $V=\oplus_{\alpha\in A}V(\alpha)$ and the fusion rule $V(\alpha)\times V(\beta)=V(\alpha+\beta)$ holds for all $\alpha,\beta\in A$.
So we often denote $V=\oplus_{W\in S_V}W$, where $S_V$ is the set of isomorphism classes of irreducible $V(0)$-modules $V(\alpha)$, $\alpha\in A$.

In this subsection, let $V(0)$ be a simple VOA satisfying the following conditions:
\begin{enumerate}[{\rm (a)}]
\item $V(0)$ has finitely many irreducible modules.
\item Any irreducible $V(0)$-module is a simple current.
\item The associativity of the fusion rules of $V(0)$ holds.
\end{enumerate}

We study a simple VOA $V$ containing $V(0)$ as a full subVOA.
Let $S_0$ be the set of all isomorphism classes of irreducible $V(0)$-modules.
Then the following proposition holds.

\bp\label{LScE} \begin{enumerate}
\item $S_0$ is a finite abelian group under the fusion rules.
\item Suppose that $V$ is a direct sum of irreducible $V(0)$-modules.
Then $V$ is a simple current extension of $V(0)$.
\end{enumerate}
\ep
\proof\ The property (b) shows that $\times$ is a binary operation on $S_0$.
Let us show that $S_0=(S_0,\times)$ is a finite abelian group.
It is obvious that the isomorphism class of $V(0)$ is the identity element.
The properties (a) and (c) show that $S_0$ is finite and associative respectively.
For an irreducible module $M$ of a VOA let $M^\prime$ be the contragraduent module of $M$ (cf. \cite{FHL}).
Then the fusion rule $M\times M^\prime =V(0)^\prime$ holds.
Moreover $V(0)^\prime\times N^\prime=V(0)$, where $N=V(0)^\prime\times V(0)^\prime$.
Hence $M\times (M^\prime\times N^\prime)=V(0)$.
This shows that any element of $S_0$ has its inverse.
By Proposition 5.4.7 in \cite{FHL}, the fusion rules are commutative, so $S_0$ is abelian.
Hence we obtain (1).

By the assumption of (2), $V=\oplus_{W\in S_0}\mu_W W$ as $V(0)$-modules, where $\mu_W$ is the multiplicity.
The simplicity of $V$ shows that $\mu_{V(0)}=1$.
Set $S_1=\{W\in S_0|\ \mu_W\neq 0\}$ and let $\tilde{S}_1$ be a subgroup of $S_0$ generated by $S_1$.
Then the abelian group $\tilde{S}_1^*$ of irreducible characters of $\tilde{S}_1$ acts faithfully on $V$ as automorphisms of a VOA.
Since $V(0)$ is a full subVOA of $V$, $V(0)$ is the fixed points of $\tilde{S}_1^*$.
By \cite{DM1} $\mu_W\in\{0,1\}$ and $\tilde{S}_1=S_1$, which shows (2).\qe

\br If $V(0)$ is rational then any VOA $V$ containing $V(0)$ as a full subVOA satisfies the assumption of (2).
\er

In the next subsection, we give examples of VOAs satisfying (a)-(c).

\subsection{$V_L^+$ for a $2$-elementary totally even lattice}
Let $L$ be a $2$-elementary totally even lattice of rank $n$, namely the dual lattice $L^*=\{\alpha\in \R\otimes_\Z L|\ \langle \alpha,L\rangle\subset\Z\}$ is a sublattice of $L/2$ and both $L$ and $\sqrt2L^*$ are even, where $\langle\cdot,\cdot\rangle$ is a positive-definite symmetric bilinear form on $\R\otimes_\Z L$.
In this subsection, we review the properties of the VOA $V_L^+$.

Let $\hat{L}$ be the central extension of $L$ by $\langle\kappa_L|\ \kappa_L^2=1\rangle$ $$1\to\langle\kappa_L\rangle\to\hat{L}\ \bar{\to}\ L\to1$$
such that $[a,b]=\kappa_L^{\langle\bar{a},\bar{b}\rangle}$ for $a,b\in\hat{L}$.
Let $\theta_L$ be the automorphism of $\hat{L}$ defined by $\theta_L(a)=a^{-1}\kappa_L^{\langle \bar{a},\bar{a}\rangle/2}$, $a\in\hat{L}$.
Set $K_L=\{a^{-1}\theta_L(a)|\ a\in\hat{L}\}$.
Then $K_L$ is normal; Consider $\hat{L}/K_L$.
Let $T$ be an irreducible $\hat{L}/K_L$-module on which $\kappa_L K_L$ acts by $-1$.
Since $L$ is $2$-elementary totally even, the center $Z(\hat{L}/K_L)$ of $\hat{L}/K_L$ is isomorphic to $2L^*/2L\times\langle\kappa_LK_L\rangle$.
By Theorem 5.5.1 in \cite{FLM} $T$ is characterized by an irreducible character of $2L^*/2L$.
For any irreducible character $\chi$ of $2L^*/2L$, there exists a unique element $\lambda+L$ in $L^*/L$ such that $\chi(\mu)=(-1)^{\langle\lambda,\mu\rangle}$ for all $\mu\in2L^*/2L$.
We denote such an irreducible character by $\chi_\lambda$ and the corresponding irreducible $\hat{L}/K_L$-module by $T_{\chi_\lambda}$.

Let $V_L$ denote the VOA associated with $L$ (cf. \cite{Bo,FLM}).
Let $\theta_{V_L}$ be an involution of $V_L$ induced by $\theta_L$.
Then $V_L^+=\{v\in V_L|\ \theta_{V_L}(v)=v\}$ is a subVOA of $V_L$.

Applying the results of \cite{DN,AD} to our case, we obtain the following proposition:
\bp\label{PIM}{\rm \cite{DN,AD}} Let $L$ be a $2$-elementary totally even lattice.
Then any irreducible $V_L^+$-module is isomorphic to one of $V_{\lambda+L}^\pm$ and $V_L^{T_{\chi_\lambda},\pm}$ $(\lambda\in L^*/L)$.
In particular, $V_L^+$ has exactly $2^{m+2}$ non-isomorphic irreducible modules, where $m=|L^*/L|$.
\ep

In this paper, we use the following notation for the isomorphism classes of irreducible $V_L^+$-modules; Let $[\lambda]^\pm$ and $[\chi_\lambda]^\pm$ denote the isomorphism classes of $V_{\lambda+L}^\pm$ and $V_{L}^{T_{\chi_\lambda},\pm}$ for $\lambda\in L^*/L$ respectively.

\bn\label{RIM}{\rm The irreducible $V_L^+$-modules $V_L^{T_{\chi_\lambda},\pm}$ are the $\pm1$-eigenspaces of the involution induced from the identity map on $T_{\chi_\lambda}$.
This notation is different from that of \cite{FLM}.}
\en

The fusion rules of $V_L^+$ are explicitly described in \cite{Ab,ADL}.
In particular, the following proposition holds.

\bp\label{PF} {\rm \cite{Ab,ADL}} Let $L$ be a $2$-elementary totally even lattice.
Then any irreducible $V_L^+$-module is a simple current.
\ep

So, we can apply the result in the previous section to $V_L^+$ by the following lemma.

\bl\label{LEx} Let $L$ be a $2$-elementary totally even lattice.
Then $V_L^+$ satisfies (a)-(c).
\el
\proof\ Proposition \ref{PIM} and \ref{PF} shows that $V_{L}^+$ satisfies (a) and (b).
By Theorem 5.18 of \cite{ADL}, the fusion rules of $V_L^+$ are associative, so $V_L^+$ satisfies (c).\qe

We denote by $S_L$ the set of all isomorphism classes of irreducible $V_L^+$-modules.
Then $S_L$ forms an elementary abelian $2$-group under the fusion rules $\times$ by \cite{Ab,ADL} and Proposition \ref{LScE}.
We often view $S_L$ as an $(m+2)$-dimensional vector space over $\F_2$, where $m=|L^*/L|$.

\bp\label{PQF}{\rm \cite{Sh2}} Suppose that the rank of $L$ is $8$ or $16$.
Then the following map $q_L$ from $S_L$ to $\F_2$ is a non-singular quadratic form on $S_L$:
\begin{eqnarray*}
q_L(W)=\left\{\begin{array}{cl}
 \mbox{$0$} & \mbox{${\rm if}\ \dim_*(W)\in\Z[[q]]$},\\
 \mbox{$1$} & \mbox{${\rm if}\ \dim_*(W)\in q^{1/2}\Z[[q]]$},
\end{array}
\right.
\end{eqnarray*}
where $\dim_*(W)=\sum_{j\in\Q}(\dim M_{j})q^j$ for a representative $M=\oplus_{j\in\Q}M_j$ of $W$.
\ep

Finally, we summarize some facts on the automorphism group of $V_L^+$.
Let $O(L)$ denote the group of all linear automorphisms of $L$ preserving the inner product $\langle\cdot,\cdot\rangle$.
For an automorphism $g$ of $\hat{L}$, let $\bar{g}$ denote the linear automorphism of $L$ defined by $\bar{g}(\bar{a})=\overline{g(a)}$, $a\in\hat{L}$.
Set $O(\hat{L})=\{g\in\Aut(\hat{L})|\ \bar{g}\in O(L)\}$.
For $\beta\in L^*/2L^*$, let $f_\beta$ denote the group homomorphism from $L$ to $\Z_2$ given by
\begin{eqnarray}
f_\beta:\gamma\mapsto \langle \beta,\gamma\rangle\mod2.\label{Def:f}
\end{eqnarray}
Then $\Hom(L,\Z_2)=\{f_\beta|\ \beta\in L^*/2L^*\}$.
We regard $f_\beta$ as an automorphism of $\hat{L}$ as follows: $$f_\beta: a\mapsto\kappa_L^{f_\beta(\bar{a})}a.$$
Hence we obtain an embedding $\Hom(L,\Z_2)\hookrightarrow O(\hat{L})$.
By \cite[Proposition 5.4.1]{FLM}, the following sequence
\begin{eqnarray} 1\to\Hom(L,\Z_2)\hookrightarrow O(\hat{L})\ \bar{\to}\ O(L)\to 1\label{Eq:Aut}\end{eqnarray}
is exact.
By \cite[Corollary 10.4.8]{FLM}, $O(\hat{L})$ acts faithfully on $V_L$.
Moreover, $O(\hat{L})/\langle\theta_{V_L}\rangle$ acts faithfully on $V_L^+$.

On the other hand, in Chapter 11 of \cite{FLM} some automorphisms of $V_L^+$ are explicitly constructed.
For Construction B, a procedure making lattices from codes, see (\ref{Eq:CB}) in Section \ref{SLL}.

\bp\label{PExA}{\rm \cite{FLM}} Let $L$ be an even lattice obtained by Construction B.
Then $V_L^+$ has automorphisms not belonging to $O(\hat{L})/\langle\theta_{V_L}\rangle$.
\ep

In \cite{Sh2,Sh3}, the automorphism group $\Aut(V_L^+)$ of $V_L^+$ is studied by using its action on $S_L$.
We collect results necessary in this paper.

\bp\label{PSh}{\rm \cite{Sh2,Sh3}} Let $L$ be a $2$-elementary totally even lattice of rank $n$ without roots.
Then the following hold:
\begin{enumerate}
\item $\Aut(V_L^+)\cong O(\hat{L})/\langle\theta_{V_L}\rangle$ if and only if $\{\lambda+L\in L^*/L|\ \#(\lambda+L)_2=2n\}=\phi$.
In particular, if $L$ is unimodular then $\Aut(V_L^+)\cong O(\hat{L})/\langle\theta_{V_L}\rangle$.\label{PSh2}
\item For $g\in\AutO(\hat{L})/\langle\theta_{V_L}\rangle$, we have 
\begin{eqnarray*}
{}\{[\lambda]^{\pm}\circ g\}&=&\{[\bar{g}^{-1}(\lambda)]^\pm\},\ \lambda\in L^*/L,\\
{}[0]^{\pm}\circ g&=&[0]^\pm.
\end{eqnarray*}
Moreover for $f_\beta\in\Hom(L,\Z_2)$
\begin{eqnarray*}
[\lambda]^{\pm}\circ f_\beta&=&
\left\{\begin{array}{cl}
 \mbox{$[\lambda]^\pm$} & \mbox{${\rm if}\ \langle\beta, \lambda\rangle\in\Z$},\\
 \mbox{$[\lambda]^\mp$} & \mbox{${\rm if}\ \langle\beta, \lambda\rangle\in\Z+1/2$},
\end{array}
\right.\\
{}[\chi_\lambda]^{\pm}\circ{f_\beta}&=&[\chi_{\lambda+\beta}]^\pm.
\end{eqnarray*}\label{PSh3}
\item $\Aut(V_{\sqrt2E_8}^+)\cong O^+(10,2)$ and $\Aut(V_{\BW}^+)\cong 2^{16}\cdot\Omega^+({10},2)$, where $\BW$ is the Barnes-Wall lattice of rank $16$.\label{PSh4}
\item Suppose that $L\cong\sqrt2E_8$ or $\BW$.
Then $S_L$ is decomposed into three-orbits under the action of $\Aut (V_L^+)$.
Moreover, $\Aut(V_{L}^+)/O_2(\Aut(V_{L}^+))$ acts faithfully on $S_L$.\label{PSh5}

\end{enumerate}
\ep

\subsection{Moonshine module}
In this subsection, we recall some facts on the moonshine module from \cite{FLM}.

Let $\Lambda$ be the Leech lattice.
Since $\Lambda$ is unimodular, $\hat{\Lambda}/K_\Lambda$ is isomorphic to the extraspecial $2$-group $2_+^{1+24}$.
Hence $\hat{\Lambda}/K_\Lambda$ has a unique faithful irreducible module $T$ on which the central element $\kappa_\Lambda K_\Lambda$ acts by $-1$.
The moonshine module $V^\natural$ is defined by $V^\natural=V_\Lambda^+\oplus V_\Lambda^{T,-}$.

\bt\label{TMO}{\rm \cite{FLM}} The $V_\Lambda^+$-module $V^\natural$ has a unique vertex operator algebra structure up to isomorphism extending its $V_\Lambda^+$-module structure.
\et

Let $z_{V^\natural}$ denote the automorphism of $V^\natural$ which acts by $1$ on $V_\Lambda^+$ and by $-1$ on $V_\Lambda^{T,-}$.
In Section 10.3 of \cite{FLM}, it was proved that $\Aut(V^\natural)$ contains a subgroup isomorphic to a non-split central extension $C$ of $O(\hat{\Lambda})/\langle\theta_{V_\Lambda}\rangle$ by $\langle z_{V^\natural}\rangle$.
By the sequence (\ref{Eq:Aut}), we have a canonical map\ $\bar{}$\ from $C$ to $O(\Lambda)/\langle-1\rangle$.
Its kernel is a central extension of $\Hom(\Lambda,\Z_2)$ by $\langle z_{V^\natural}\rangle$, which is isomorphic to $\hat{\Lambda}/K_\Lambda$.
So, the sequence of groups
\begin{eqnarray}
1\to\hat{\Lambda}/K_\Lambda\hookrightarrow C\ \bar{\to}\ \AutO(\Lambda)/\langle-1\rangle\to1\label{Def:C}
\end{eqnarray}
is exact.

\bn We will show that $C$ is the centralizer of $z_{V^\natural}$ in $\Aut(V^\natural)$ later without properties of the Monster.
\en

In Chapter 11 and 12 in \cite{FLM}, an automorphism $\sigma$ of $V^\natural$ not belonging to $C$ was explicitly constructed.
By the fact that the degree $2$ subspace of $V^\natural$ has the algebraic structure given by \cite{Gr1}, the following theorem was shown.

\bt\label{MTFLM}{\rm \cite{FLM}} The automorphism group of $V^\natural$ is generated by $C$ and $\sigma$, and it is isomorphic to the Monster.
\et

\bn In this paper we study symmetries of $V^\natural$ without Theorem \ref{MTFLM}.
\en

\subsection{Sublattices of the Leech lattice}\label{SLL}
In this subsection, we study some sublattices of the Leech lattice and their automorphism groups.
For the precise definitions of the Golay code and the Leech lattice, see \cite{CS}.

We set $\Omega_n=\{1,2,\dots,n\}$.
We often identify the power set $\mathcal{P}(\Omega_n)$ of $\Omega_{n}$ with an $n$-dimensional vector space over $\F_2$.
Let $(\cdot,\cdot)$ denote the natural inner product on $\mathcal{P}(\Omega_n)$.
Let $\{\alpha_j\mid j\in\Omega_n\}$ be an orthogonal basis of $\R^{n}$ of norm $2$.
For a linear binary code $\mathcal{C}$ of length $n$, the lattice
\begin{eqnarray}
L_B(\mathcal{C})=\sum_{c\in \mathcal{C}}\Z\frac{1}{2}\alpha_c+\sum_{j,k\in\Omega_n}\Z(\alpha_j+\alpha_k)\label{Eq:CB}
\end{eqnarray}
is said to be {\it the lattice obtained by Construction B from $\mathcal{C}$}, where $\alpha_c=\sum_{j\in c}\alpha_j$.

Let $G_{24}$ be the extended Golay code, which is the unique $12$-dimensional doubly even binary code of length $24$ with minimum weight $8$.
A codeword of $G_{24}$ of weight $8$ is called an {\it octad}.
A partition $\{P_l\mid l\in\Omega_6\}$ of $\Omega$ is called a {\it sextet} if $|P_l|=4$ and $P_l\cup P_k$ is an octad for all $l\neq k$.

Let $\Lambda$ be the Leech lattice, which is the unique positive-definite even unimodular lattice of rank $24$ without roots up to isomorphism.
The automorphism group of $\Lambda$ is denoted by $Co_0$ and its quotient by $\langle-1\rangle$ is denoted by $Co_1$.
We set $\Omega=\Omega_{24}$.
Then the lattice
\begin{eqnarray*}
&&L_B(G_{24})+\Z(\frac{\alpha_\Omega}{4}-\alpha_1)\notag\\
&=&\sum_{c\in G_{24}}\Z\frac{\alpha_c}{2}+\sum_{j,k\in\Omega}\Z(\alpha_j+\alpha_k)+\Z(\frac{\alpha_\Omega}{4}-\alpha_1).
\end{eqnarray*}
is an even unimodular lattice of rank $24$ without roots, which is the Leech lattice $\Lambda$.
We use this expression of the Leech lattice in this paper.

Let us consider $4$ full sublattices $\Lambda(i)$ ($i=1,2,3,5$) of $\Lambda$.
We set $\Lambda(1)=\Lambda$ and $\Lambda(2)=L_B(G_{24})$.
We fix a sextet $\{P_l^1|\ l\in\Omega_{6}\}$.
Set $C(3)=\{c\in G_{24}|\ (P^1_1,c)=0\}$ and $\Lambda(3)=L_B(C(3))$.
We denote $O_s=P_{2s}^1\cup P_{2s-1}^1$ for $s=1,2,3$.
Let $\{P_l^j|\ l\in\Omega_6\}$ ($j=2,3$) be distinct sextets such that $P_{2s-1}^j, P_{2s}^j\subset O_s$ and $\F_2\langle P_{2s}^j,O_s|\ j=1,2,3\rangle\subset\mathcal{P}(O_s)$ is isomorphic to the extended Hamming code of length $8$ for $s=1,2,3$.
Let $C(5)$ be the subcode of $C(3)$ defined by $C(5)=\{c\in G_{24}|\ (P^j_1,c)=0\ {\rm for}\ j=1,2,3\}$.
We set $\Lambda(5)=L_B(C(5))$.
Then it is easy to check that $$|\Lambda/\Lambda(i)|=2^{i-1}.$$

We now consider other expressions of the lattices $\Lambda(i)$.
We set $D_1=2\Lambda$, $D_2=\Z\langle 2\alpha_1,D_1\rangle$, $D_3=\Z\langle \alpha_{P_{1}^1},D_2\rangle$ and $D_5=\Z\langle \alpha_{P_1^j}, D_3|\ j=2,3\rangle$.
Then
\begin{eqnarray}
\Lambda(i)=\{v\in\Lambda|\ \langle v,D_i\rangle\subset2\Z\}=D_i^*/2.\label{Eq:Lambda235}
\end{eqnarray}

Recall that an even lattice $L$ is $2$-elementary if the dual lattice $L^*=\{\alpha\in \R\otimes_\Z L|\ \langle \alpha,L\rangle\subset\Z\}$ is a sublattice of $L/2$, and is totally even if $\sqrt2L^*$ is even.

\bp\label{PL1} The lattice $\Lambda(i)$ is $2$-elementary totally even for $i=1,2,3,5$.
\ep
\proof\ By the inclusion $\Lambda(1)\supset\Lambda(2)\supset\Lambda(3)\supset\Lambda(5)$ and the unimodular property of $\Lambda(1)$, it suffices to consider $\Lambda(5)$.
It is easy to see that
\begin{eqnarray}
\Lambda(5)^*=\Z\langle\Lambda,\alpha_1,\frac{\alpha_{P_1^j}}{2}|\ j=1,2,3\rangle.\label{Eq:GenLa5}
\end{eqnarray}
Hence $\Lambda(5)^*/\Lambda(5)\cong 2^{8}$.
In particular $2\Lambda(5)^*\subset\Lambda(5)$.
It is easy to see that the norms of the generators of $\Lambda(5)^*$ in (\ref{Eq:GenLa5}) are integer and that their inner products are in $\Z/2$.
Hence the norm of each vector of $\Lambda(5)^*$ is integer.
Therefore $\Lambda(5)$ is $2$-elementary totally even.
\qe

Let $O(L)$ denote the group of all linear automorphisms of $L$ preserving the inner product $\langle\cdot,\cdot\rangle$.
For a sublattice $N$ of a lattice $L$, we denote by $O(L;N)$ the subgroup of $O(L)$ consisting of all automorphisms preserving $N$.
Let us study $O(\Lambda;\Lambda(i))$.
By (\ref{Eq:Lambda235}), $O(\Lambda;\Lambda(i))=\{g\in O(\Lambda)|\ g(D_i)=D_i\}$.
Let $\rho$ be the canonical homomorphism from $\Lambda$ to $\Lambda/2\Lambda$.
Note that $Co_1$ acts faithfully on $\Lambda/2\Lambda$.
In \cite{ATLAS,Wi} the stabilizer of $\rho(D_i)$ in $Co_1$ is described.
In particular, we obtain the following lemma.

\bl\label{LL2} {\rm \cite{ATLAS,Wi}} For $i=2,3,5$, the group $\AutO(\Lambda;\Lambda(i))$ is a maximal subgroup of $\AutO(\Lambda)$ and its shape is given as follows:
\begin{eqnarray*}
\AutO(\Lambda;\Lambda(2))&\cong& 2^{12}:M_{24},\\
\AutO(\Lambda;\Lambda(3))&\cong& 2^{5+12}\cdot(L_2(2)\times 3Sym_6),\\
 \AutO(\Lambda;\Lambda(5))&\cong& 2^{3+12}:(L_4(2)\times Sym_3).
\end{eqnarray*}
\el

Let $t_i$ be the canonical homomorphism $O(\Lambda;\Lambda(i))\to {\rm GL}(\Lambda/\Lambda(i))\cong L_{i-1}(2)$.
By the lemma above, we obtain the following lemma.

\bl\label{PL3} For each $i$, $t_i$ is surjective and its kernel $\Ker t_i$ is given as follows:
\begin{eqnarray*} 
\Ker t_2\cong2^{12}:M_{24},\ \Ker t_3\cong2^{5+12}.3Sym_6,\ \Ker t_5\cong 2^{3+12}:Sym_3.
\end{eqnarray*}
\el

We consider other full sublattices of $\Lambda$.
We recall fundamental automorphisms of $\Lambda$.
For $c\in\Omega$, let $\varepsilon_c$ denote the linear automorphism of $\R^{24}$ defined by setting
\begin{eqnarray*}
\varepsilon_c(\alpha_j)=\left\{\begin{array}{cc}
 \mbox{$-\alpha_j$} & \mbox{${\rm if}\ j\in c$},\\
 \mbox{$\alpha_j$} & \mbox{${\rm if}\ j\notin c$}.
 \end{array}
\right.
\end{eqnarray*}
It is easy to see that $\varepsilon_c$ is an automorphism of $\Lambda$ if and only if $c$ is in $G_{24}$.
Fix an octad $c$ of $G_{24}$ and set
\begin{eqnarray}
U^1=\{v\in\Lambda|\ \varepsilon_{c}(v)=-v\},\ U^2=\{v\in\Lambda|\ \varepsilon_{c}(v)=v\}.\label{Eq:DefU^i}
\end{eqnarray}
Then $U=U^1\oplus U^2$ is a sublattice of $\Lambda$.
It is well known that $U^1\cong\sqrt2E_8$ and $U^2\cong \Lambda_{16}$, where $\Lambda_{16}$ is the Barnes-Wall lattice of rank $16$.
It is easy to check that $\sqrt2E_8$ and $\BW$ are $2$-elementary totally even.
The determinant of $U$ is $2^{16}$ since that of $U^i$ is $2^8$ for $i=1,2$.
Thus $|\Lambda/U|=2^8$ and $|U/2\Lambda|=2^{16}$.
Furthermore the automorphism group of $U$ is described as follows:

\bl\label{LL4} The automorphism group $\AutO(U)$ of $U$ is isomorphic to the direct product of the groups $\AutO(U^1)$ and $\AutO(U^2)$.
\el
\proof\ We regard $\AutO(U^i)$ as a subgroup of $\AutO(U)$.
Then $\AutO(U)\supseteq\AutO(U^1)\times\AutO(U^2)$.
Hence it suffices to show that $\AutO(U)$ preserves both $U^1$ and $U^2$.
Since $\sqrt2E_8$ is indecomposable and $\AutO(U)$ preserves the inner product, we have $g(U^1)=U^1$ or $g(U^1)\subset U^2$.
If $g(U^1)\subset U^2$ then there exists a sublattice $L$ of $U^2$ such that $U^2=g(U^1)\oplus g(L)$, which contradicts the indecomposability of $\BW$.
Therefore we have $g(U^1)=U^1$ and $g(U^2)=U^2$.\qe

\section{Isomorphisms between extensions of a VOA graded by a finite abelian group}

Let $A$ be a finite abelian group and let $V$ be a simple $A$-graded VOA, namely $V=\oplus_{\alpha\in A} V(\alpha)$ and $u_{n}v\in V(\alpha+\beta)$ for any $u\in V(\alpha)$, $v\in V(\beta)$ and $n\in\Z$.
Then $V(\alpha)$ is a $V(0)$-module for all $\alpha\in A$.
In this article, we always assume that $V(\alpha)\neq0$ for all $\alpha\in A$.
So the group $A^*$ of all irreducible characters of $A$ acts faithfully on $V$ as automorphisms: for $\chi\in A^*$,  $\chi (v)=\chi(\alpha) v$, $v\in V(\alpha)$.
Clearly $V(\alpha)$ is an eigenspace of $A^*$ for all $\alpha\in A$ and $V(0)$ is the fixed points of $A^*$.
By \cite{DM1}, $V(\alpha)$, ($\alpha\in A$) are non-isomorphic and irreducible.
The following theorem is a slight generalization of Proposition 5.3 in \cite{DM2}.

\bt\label{PDM2} Let $V=\oplus_{\alpha\in A}V(\alpha)$ and $V^\prime=\oplus_{\alpha\in A}V^\prime(\alpha)$ be simple VOAs graded by a finite abelian group $A$.
Suppose that $V(0)= V^\prime(0)$ and that the fusion rule $V(\alpha)\times V(\beta)=V(\alpha+\beta)$ holds for all $\alpha,\beta\in A$.
Let $g$ be an automorphism of $V(0)$ which maps the set of isomorphism classes of $\{V(\alpha)|\ \alpha\in A\}$ to those of $\{V^\prime(\alpha)|\ \alpha\in A\}$.
Then there exists an isomorphism $\tilde{g}$ from $V^\prime$ to $V$ such that $\tilde{g}_{|V(0)}=g$.
\et
\proof\ For $\alpha\in A$, let $\varphi_{\alpha}$ be an isomorphism from $W=V(\alpha)\circ g$ to $V(\alpha)$ such that $$ \varphi_{\alpha}Y_{W}(v,z)=Y_{V(\alpha)}(gv,z)\varphi_{\alpha}$$ for $v\in V(0)$.
We regard $\varphi_{\alpha}$ as a linear isomorphism from $V^\prime(\beta)$ to $V(\alpha)$, where $V^\prime(\beta)\cong V(\alpha)\circ g$.
In particular we take $\varphi_{V(0)}$ as $g$.
Then we obtain an isomorphism $\varphi_{V^\prime}=\oplus_{\alpha\in A}\varphi_{\alpha}$ of $V^0$-modules from $V^\prime$ to $V$.
Set $\tilde{Y}_{V}(v,z)=\varphi_{V^\prime}Y_{V^\prime}(v,z)\varphi_{V^\prime}^{-1}$, $v\in V(0)$.
Then $\varphi_{V^\prime}$ is an isomorphism of VOAs from $(V^\prime,Y_{V^\prime})$ to $(V,\tilde{Y}_V)$.

Clearly $(V,\tilde{Y}_V)$ is a $A$-graded VOA: $V=\oplus_{\alpha\in A} \tilde{V}(\alpha)$.
Moreover $V(\alpha)$ is isomorphic to $\tilde{V}(\alpha)$ as $V(0)$-modules for all $\alpha\in A$.
By Proposition 5.3 in \cite{DM2}, there exists an isomorphism $\psi$ of VOAs from $(V,\tilde{Y}_V)$ to $(V,Y_V)$ such that $\psi_{|V(0)}$ is the identity map.
Therefore we obtain an isomorphism $\tilde{g}=\psi\circ \varphi$ of VOAs from $V^\prime$ to $V$ such that $\tilde{g}_{|V(0)}=g$.\qe

Let $S_A$ be the set of the isomorphism classes of the irreducible $V(0)$-modules $V(\alpha)$, $(\alpha\in A)$.
For an automorphism $g$ of $V(0)$, we set $S_A\circ g=\{W\circ g|\ W\in S_A\}$.
Then we obtain the restriction homomorphisms
\begin{eqnarray*}
\Phi^N_A&:N_{\Aut(V)}(A^*)&\to H^N_A,\\ \Phi^C_A&:C_{\Aut(V)}(A^*)&\to H^C_A,
\end{eqnarray*}
where 
\begin{eqnarray*}
H^N_A&=&\{h\in\Aut(V(0))|\ S_A\circ h=S_A\},\\
H^C_A&=&\{h\in\Aut(V(0))|\ W\circ h=W\ {\rm for\ all}\ W\in S_A\}.
\end{eqnarray*}

Applying Theorem \ref{PDM2} to the case $V=V^\prime$, we show that $\Phi^N_A$ is surjective.
Since each $V(\alpha)$ is irreducible, $\Ker\Phi^N_A=A^*$.
By similar arguments, $\Phi^C_A$ is surjective, and $\Ker\Phi^C_A=A^*$.
Hence we obtain the following corollary.

\bc\label{PSC}\label{CSC}{\rm (cf. \cite[Theorem 3.3]{Sh2})} Let $V=\oplus_{\alpha\in A}V(\alpha)$ be a simple VOA graded by a finite abelian group $A$.
Suppose that the fusion rule $V(\alpha)\times V(\beta)=V(\alpha+\beta)$ holds for all $\alpha,\beta\in A^*$.
Then the restriction homomorphism $\Phi^N_A$ and $\Phi^C_A$ are surjective and $\Ker\Phi^C_A=\Ker\Phi^N_A=A^*$.
\ec

By the corollary above, we obtain the following exact sequences:
\begin{eqnarray}0\to A^*\to N_{\Aut(V)}(A^*)\to H_A^N\to 1,\label{Eq:NA}\\ 0\to A^*\to C_{\Aut(V)}(A^*)\to H_A^C\to 1.\label{CA}\end{eqnarray}
This shows that the normalizer and centralizer of $A^*$ in $\Aut(V)$ are described in terms of $\Aut(V(0))$.

\section{Simple current extensions of $V_L^+$}
Let $L$ be the even lattice isomorphic to either $\sqrt2E_8$ or $\BW$.
In this section, we classify VOAs containing $V_{L}^+$ as a full subVOA by using the symmetries of $V_{L}^+$.

Since $L$ is $2$-elementary totally even, $V_L^+$ satisfies the conditions (a)-(c) in Section \ref{Ss} by Lemma \ref{LEx}.
Moreover, $L$ has a $4$-frame, an orthogonal basis of $\R\otimes_\Z L$ of norm $4$, so $V_{L}^+$ is a framed VOA.
In particular $V_{L}^+$ is rational by Theorem 2.12 in \cite{DGH}.
Let $S_L$ be the set of all isomorphism classes of irreducible $V_L^+$-modules.
Then $S_L$ is an elementary abelian $2$-group of order $2^{10}$ by Lemma \ref{LScE} (1).
Moreover $S_L$ has a quadratic form $q$ given in Proposition \ref{PQF}.
Let $V$ be a simple VOA containing $V_L^+$ as a full subVOA.
Then by Proposition \ref{LScE} (2), $V$ is a simple current extension of $V_L^+$: $V=\oplus_{W\in S_V}W$, where $S_V$ is a subgroup of $S_L$.
Clearly degrees of each element of $S_V$ belong to the set of integers.
Hence we obtain the following lemma.
\bl Let $V=\oplus_{W\in S_V}W$ be a simple current extension of $V_L^+$.
Then $S_V$ is a totally singular subspace of $S_L$.
\el

Now we recall fundamental results on orthogonal groups (cf. \cite[Theorem 9.4.3, Theorem 9.4.8]{BCN}).
\bp\label{PGT} Let $V$ be a $10$-dimensional vector space over $\F_2$ with a non-singular quadratic form of plus type.
Let $O(V)$ be the orthogonal group of degree $10$ on $V$.
For $0\le i\le 5$, let $\mathcal{S}_i$ denote the set of all $i$-dimensional totally singular subspaces of $V$.
Then the following hold:
\begin{enumerate}
\item For $0\le i\le  5$, $O(V)\cong O^+(10,2)$ acts transitively on $\mathcal{S}_i$ .
\item For $0\le i\le 4$, the commutator subgroup $O(V)^\prime\cong \Omega^+(10,2)$ acts transitively on $\mathcal{S}_i$
\item $\mathcal{S}_5$ is decomposed into $2$-orbits under the action of $O(V)^\prime$.
\end{enumerate}
\ep

By Proposition \ref{PSh} (3), Theorem \ref{PDM2} and Proposition \ref{PGT}, the number of isomorphism classes of simple current extensions of $V(0)$ is less than or equal to $6$ if $L=\sqrt2E_8$, and $7$ if $L=\BW$.

Let $N$ be an even overlattice of $L$ of index $2^j$ $(0\le j\le 4)$.
Then $V_N$ is a simple current extension of $V_L^+$ graded by an elementary abelian $2$-group of order $2^{j+1}$.
We note that $N$ is unique up to isomorphism for each $j$ except that $E_8\oplus E_8$ and $\Gamma_{16}$ are non-isomorphic even overlattices of $\BW$ of index $2^4$.
Therefore we obtain the following theorem.
\bt Let $L$ be a lattice isomorphic to either $\sqrt2E_8$ or $\BW$.
\begin{enumerate}
\item Any simple VOA containing $V_L^+$ as a proper full subVOA is isomorphic to $V_N$, where $N$ is an even overlattice of $L$.
In particular the number of isomorphism classes of simple current extensions of $V_L^+$ is $6$ if $L\cong \sqrt2E_8$, and $7$ if $L\cong \BW$.
\item Let $S_V$ be a subset of $S_L$.
The $V_L^+$-module $V=\oplus_{W\in S_V}\mu_WW$ has a VOA structure if and only if $S_V$ is a totally singular subspace of $S_L$ and $\mu_W=1$ for all $W\in S_V$.
In particular the VOA structure on $V$ is unique up to isomorphism.
\end{enumerate}
\et

\section{Normalizers of elementary abelian $2$-groups in the automorphism group of the moonshine module}\label{SNSOM}
In this section, we consider some subVOAs $V_L^+$ of the moonshine module $V^\natural$ and describe the normalizers in $G=\Aut(V^\natural)$ of associated elementary abelian $2$-groups in terms of $\Aut(V_L^+)$.
Moreover, we show that our normalizers form a Monster amalgam.

\subsection{Decomposition of $V^\natural$ as irreducible $V_L^+$-modules}
Let $L$ be a $2$-elementary totally even full sublattice of the Leech lattice $\Lambda$.
In this subsection, we decompose the moonshine module $V^\natural$ into irreducible $V_L^+$-modules.
Since VOA $V_L^+$ is a subVOA of $V_\Lambda^+$, it suffices to decompose $V_\Lambda^+$ and $V_\Lambda^{T,-}$.

The following proposition is easy.
\bp\label{LD1} Let $L$ be a $2$-elementary totally even full sublattice of the Leech lattice $\Lambda$.
Then the VOA $V_\Lambda^+$ decomposes into irreducible $V_L^+$-modules as follows:
\begin{eqnarray*}
V_\Lambda^+\cong\bigoplus_{\lambda+L\in\Lambda/L}V_{\lambda+L}^+.
\end{eqnarray*}
\ep

Now we decompose $V_\Lambda^{T,-}$ into $V_L^+$-modules as in \cite{Sh1}.

\bp\label{LD2} Let $L$ be a $2$-elementary totally even full sublattice of the Leech lattice $\Lambda$.
Then the $V_\Lambda^+$-module $V_\Lambda^{T,-}$ decomposes into irreducible $V_L^+$-modules as follows:
\begin{eqnarray*}
V_\Lambda^{T,-}\cong\bigoplus_{\lambda+L\in \Lambda/L}V_{L}^{T_{\chi_{\lambda}},-}.
\end{eqnarray*}
\ep
\proof\ In order to decompose $V_\Lambda^{T,-}$, we will decompose $\hat{\Lambda}/K_\Lambda$-module $T$ into $\hat{L}/K_L$-modules.
Set $N=2L^*$.
Then the center $Z(\hat{L}/K_\Lambda)$ of $\hat{L}/K_\Lambda$ is $\hat{N}/K_\Lambda$.
We choose a maximal abelian subgroup $A$ of $\hat{\Lambda}/K_\Lambda$ satisfying $\hat{N}/K_\Lambda\subset A\subset \hat{L}/K_\Lambda$.
Then 
\begin{eqnarray}
T=\oplus_{\chi\in X(A)} \C_\chi\label{Eq:DeA}
\end{eqnarray}
as irreducible $A$-modules, where $\C_\chi$ is a one-dimensional $A$-module with character $\chi$ and $X(A)$ is the set of all characters $\chi$ of $A$ with $\chi(\kappa_\Lambda K_\Lambda)=-1$.
Then $T=\oplus_{\chi\in X(\hat{L}/K_\Lambda)} m_\chi T_\chi$ as $\hat{L}/K_\Lambda$-modules, where $T_\chi$ is the irreducible $\hat{L}/K_\Lambda$-module with character $\chi$ and $m_\chi$ is its multiplicity.
Since the components in the decomposition (\ref{Eq:DeA}) are distinct $A$-modules, we have $m_\chi\in\{0,1\}$.
Let $D$ be a complement to $\langle\kappa K_\Lambda\rangle$ in $\hat{N}/K_\Lambda$.
Then $\hat{L}/D$ is isomorphic to the extraspecial $2$-group of shape $2_+^{1+2l}$, where $2^{2l}=|L/2L^*|$.
We note that a faithful irreducible module of an extraspecial $2$-group of shape $2_+^{1+2l}$ is unique up to isomorphism and that its dimension is $2^{l}$.
Then $\dim T_\chi=2^l$ for $\chi\in X(\hat{L}/K_\Lambda)$.
Moreover each character in $X(\hat{L}/K_\Lambda)$ corresponds to one of $N/2\Lambda$ by $\bar{D}= N/2\Lambda$.
Since $|N/2\Lambda|=2^{12-l}$ and $\dim T=2^{12}$, we have $m_\chi=1$ for all $\chi\in X(\hat{L}/K_\Lambda)$ by comparing the dimensions.
For a character $\chi$ of $N/2\Lambda$, there exists a unique element $\lambda$ of $\Lambda/L$ such that $\chi(\cdot)=(-1)^{\langle\lambda,\cdot\rangle}=\chi_\lambda(\cdot)$.
Since $\hat{N}/K_L\supset K_\Lambda/K_L$, we may view $T_\chi$ as a $\hat{L}/K_L$-module.
Hence we obtain the following decomposition of $T$ as a $\hat{L}/K_L$-module:
\begin{eqnarray*}
T\cong\bigoplus_{\lambda+L\in\Lambda/L}T_{\chi_{\lambda}},
\end{eqnarray*}
where $T_{\chi_{\lambda}}$ is the irreducible $\hat{L}/K_L$-module with central character $\chi_{\lambda}$.
Therefore we obtain the desired decomposition of $V_\Lambda^{T,-}$.
\qe

By the propositions above, Proposition \ref{LScE} and Lemma \ref{LEx}, we obtain the following corollary:

\bc Let $L$ be a $2$-elementary totally even full sublattice of the Leech lattice.
Then $V^\natural$ is a simple current extension of $V_L^+$ graded by an abelian group of order $2^{m+1}$, where $|L^*/L|=2^{2m}$.
\ec

\subsection{Symmetries of the moonshine module associated with $V_L^+$ (I)}
Let $\Lambda(i)$ be the $2$-elementary totally even full sublattices of the Leech lattice $\Lambda$ given in Section  \ref{SLL} for $i=1,2,3,5$.
In this subsection, we consider the abelian automorphism group of $V^\natural$ associated with the fusion rules of $V_{\Lambda(i)}^+$, and describe its normalizer in the automorphism group $G$ of the moonshine module $V^\natural$.

Let $S_{i}$ be the set of isomorphism classes of irreducible $V_{\Lambda(i)}^+$-modules which appear in the decomposition of $V^\natural$ in Proposition \ref{LD1} and \ref{LD2}.
Then $S_{i}$ forms an elementary abelian $2$-group under the fusion rules of order $2^i$ and its dual $S_i^*$ acts faithfully on $V^\natural$.

Let us study the centralizer $C_{G}(S_i^*)$ of $S_i^*$ in $G$.
We start by proving the following lemmas.

\bl The centralizer $C_{G}(S_1^*)$ of $S_1^*$ is equal to $C$.
\el
\proof\ By Lemma \ref{PSh} $\Aut(V_{\Lambda}^+)=O(\hat{\Lambda})/\langle\theta_{V_\Lambda}\rangle$.
Since the graded dimensions of $V_{\Lambda}^-$ and $V_{\Lambda}^{T,\pm}$ are distinct, $H_{S_1}^C=\{g\in\Aut(V_{\Lambda}^+)|\ V_{\Lambda}^-\circ g\cong V_{\Lambda}^-\}=\Aut(V_{\Lambda}^+)$.
By Corollary \ref{CSC}, (\ref{Eq:Aut}) and (\ref{Def:C}) we obtain this lemma.
\qe

\bl \label{LCG}For any $i$, the centralizers of $S_i^*$ in $C$ and $G$ are the same: $C_C(S_i^*)=C_{G}(S_i^*)$.
\el

\proof\ Let $g\in C_{G}(S_i^*)$.
Then $g$ commutes with $z_{V^\natural}$.
Since $S_1^*=\langle z^\natural\rangle$ and $C=C_{G}(S_1^*)$, we obtain $g\in C$.
The converse is obvious.
\qe

Let us describe the centralizer of $S_i^*$ in $C$.
Recall the canonical homomorphism $t_i: O(\Lambda;\Lambda(i))\to GL(\Lambda/\Lambda(i))$.
The image of $C_C(S_i^*)$ under the canonical map $C\to\AutO(\Lambda)$ is contained in $\Ker t_i$ since $S_i^*$ preserves  $V_{\lambda+\Lambda(i)}^\varepsilon$ for each $\lambda\in\Lambda/\Lambda(i)$, $\varepsilon\in\{\pm\}$ by Proposition \ref{PSh} (2).
Since $S_i^*$ is contained in the kernel of the homomorphism $\bar{}$, we obtain the homomorphism $C_{C}(S_i^*)/S_i^*\ \bar{\to}\ \Ker t_i$.
Let us show that this homomorphism is surjective.
For any element $g\in\Ker t_i$, there exists $h_0\in  C$ such that  $\bar{h}_0=g$.
Since $h_0$ preserves $V_{\Lambda(i)}^+$ and the decompositions in Proposition \ref{LD1} and \ref{LD2} are multiplicity free, $\{h_0(V_{\lambda+\Lambda(i)}^\varepsilon)\}=\{V_{\lambda+\Lambda(i)}^\varepsilon\}$.
By Proposition \ref{PSh} (2) and the fusion rules of $V_L^+$, we can take $h_1\in\Hom(\Lambda,\Z_2)$ such that $h_1\circ h_0(V_{\lambda+\Lambda(i)}^\varepsilon)=V_{\lambda+\Lambda(i)}^\varepsilon$ for $\varepsilon\in\{\pm\}$.
Let $\lambda$ be an element of $L^*/L$ such that $h_1\circ h_0(V_{\Lambda(i)}^{T_{\chi_0},\varepsilon})=V_{\Lambda(i)}^{T_{\chi_\lambda},\varepsilon}$.
Then $h_1\circ h_0(V_{\Lambda(i)}^{T_{\chi_\mu},\varepsilon})=V_{\Lambda(i)}^{T_{\chi_{\mu+\lambda}},\varepsilon}$ for all $\mu\in \Lambda/L$.
By Proposition \ref{PSh} (2) there exists $h_2$ such that $h_2\circ h_1\circ h_0$ preserves all $V_{\Lambda(i)}^{T,\chi_\mu}$.
Hence $h=h_2\circ h_1\circ h_0\in C_C(S_i^*)$ and $\bar{h}=g$, and the homomorphism is surjective.

Let $f_\beta\in\Hom(\Lambda,\Z_2)\subset O(\hat{\Lambda})$.
Set $U_i^*=S_i^*/\langle z_{V^\natural}\rangle$.
Then $f_\beta$ commutes $U_i^*$ if and only if $\beta$ belongs to $\Lambda(i)$ by Proposition \ref{PSh} (2), namely $f_\beta(2\Lambda(i)^*)=0$.
Then we obtain the sequence of groups
\begin{eqnarray} 1\to F_i/U_i^*\to C_{C}(S_i^*)/S_i^*\ \bar{\to}\ \Ker t_i\to 1,\label{Eq:ExC}\end{eqnarray}where $F_i=\{f\in\Hom(\Lambda,\Z_2)|\ f(2\Lambda(i)^*)=0\}$.
Since the sequence (\ref{Def:C}) is exact, the kernel of the homomorphism $\ \bar{}\ $ is $F_i/U_i^*$.
Thus the sequence (\ref{Eq:ExC}) is also exact.
By Lemma \ref{PL3} we obtain the following proposition.

\bp \label{PCG}For $i=1,2,3,5$, the sequence (\ref{Eq:ExC}) is exact.
In particular, the shapes of the centralizers are described as follows:
\begin{eqnarray*}
C_G({S_1^*})\cong&2.2^{24}.Co_1,\ &C_G({S_2^*})\cong 2^22^{22}.2^{11}.M_{24},\\
C_G({S_3^*})\cong&2^3. 2^{20}.2^{4+12}.3Sym_6,\ &C_G(S_5^*)\cong2^5.2^{16}.2^{2+12}.Sym_3.\end{eqnarray*}
\ep

Next, we discuss the normalizer $N_G(S_i^*)$ of $S_i^*$ in $G$.
Let us consider $N_C(S_i^*)$.
Since $V_{\Lambda(i)}^+$ is the set of fixed points of $S_i^*$, the image of $\ \bar{}\ $ of $N_C(S_i^*)$ is contained in $O(\Lambda;\Lambda(i))$.
By Proposition \ref{PSh} (2) $\Hom(\Lambda,\Z_2)$ preserves the set of components in Proposition \ref{LD1} and \ref{LD2}.
Hence $\Hom(\Lambda,\Z_2)\subset N_C(S_i^*)$.
Thus we obtain the sequence \begin{eqnarray}1\to \Hom(\Lambda,\Z_2)/U_i^*\to N_{C}(S_i^*)/S_i^*\ \bar{\to}\ O(\Lambda;\Lambda(i))\to 1.\label{Eq:SeqNC}\end{eqnarray}

\bl\label{LpT}
\begin{enumerate}
\item The sequence (\ref{Eq:SeqNC}) is exact.
\item The shape of $N_C(S_i^*)/C_C(S_i^*)$ is $2^{i-1}.L_{i-1}(2)$.
\end{enumerate}
\el

\proof\ By similar arguments on the case of $C_C(S_i^*)$, the group homomorphism $\bar{}$ in (\ref{Eq:SeqNC}) is surjective.
Since (\ref{Def:C}) is exact, the kernel of the homomorphism $\ \bar{}\ $ in (\ref{Eq:SeqNC}) is $\Hom(\Lambda,\Z_2)/U_i^*$.
Hence the sequence (\ref{Eq:SeqNC}) is exact.

It is easy to see that $\Hom(\Lambda,\Z_2)/F_i\cong 2^{i-1}$.
By Proposition \ref{LL2}, the canonical homomorphism $t_i:O(\Lambda;\Lambda(i))\to GL(\Lambda/\Lambda(i))\cong L_{i-1}(2)$ is surjective.
Comparing the sequences (\ref{Eq:ExC}) and (\ref{Eq:SeqNC}), we obtain $N_C(S_i^*)/C_C(S_i^*)\cong 2^{i-1}.L_{i-1}(2)$.
This completes (2).
\qe

\bt\label{T1} For $i=1,2,3,5$, the sequence 
$$ 1\to C_{G}(S_i^*)\to N_{G}(S_i^*)\to L_i(2)\to 1$$
is exact.
In particular, the shapes of $N_G(S_i^*)$ are described as follows:
\begin{eqnarray*}
N_G(S_1^*)\cong &2^{1+24}.Co_1,\ &N_G(S_2^*)\cong(2^2.2^{22}.2^{11}.M_{24}).L_2(2),\\ N_G(S_3^*)\cong&(2^3.2^{20}.2^{4+12}.3Sym_6).L_3(2),\ &N_G(S_5^*)\cong(2^5.2^{16}.2^{2+12}.Sym_3).L_5(2).
\end{eqnarray*}
\et
\proof\ The assertion is clear when $i=1$.
So we assume that $i\neq 1$.
By Lemma \ref{LCG}, Proposition \ref{PCG} and Lemma \ref{LpT} (2), we obtain a subgroup of $N_G(S_i^*)/C_G(S_i^*)\cong L_i(2)$ of shape $2^{i-1}.L_{i-1}(2)$.
We note that this subgroup is maximal in $GL(S_i^*)\cong L_i(2)$.
On the other hand, $N_G(S_i^*)$ has automorphisms not in $C_G(S_i^*)$ by Proposition \ref{PExA}.
Hence $N_G(S_i^*)/C_G(S_i^*)\cong L_i(2)$.\qe

\bn In the group theory, it was known that $S_i^*$ is $2B$-pure in the Monster for $i=1,2,3,5$ and that $N_G(S_i^*)$ is a maximal $2$-local subgroup of the Monster.
\en

\subsection{Symmetries of the moonshine module associated with $V_L^+$ (II)}\label{SDMM}
In this subsection, we consider a subVOA $V$ of $V^\natural$ isomorphic to $V_{\sqrt2E_8}^+\otimes V_{\BW}^+$.
We decompose $V^\natural$ into irreducible $V$-modules and consider the automorphism group of $V^\natural$ of shape $2^{10}$.
Then we describe the normalizer in $G=\Aut(V^\natural)$ of this elementary abelian $2$-group.
Throughout this subsection, $k$ is an element of $\{1,2\}$.

Recall the sublattice $U=U^1\oplus U^2$ of the Leech lattice $\Lambda$ associated with an octad $c$, where $U^1\cong\sqrt2E_8$ and $U^2\cong\BW$ (see (\ref{Eq:DefU^i})).
Then $U^*=(U^1)^*\oplus (U^2)^*$.
For a vector $\lambda$ of $U^*$, we denote $\lambda=\lambda_1+\lambda_2$, $\lambda_k\in (U^k)^*$.
It is easy to see that the fixed point subVOA $V$ of $V_U^+$ with respect to a lift of the automorphism $\varepsilon_{c}$ of $\Lambda$ is $V_{U^1}^+\otimes V_{U^2}^+$.

Let us study the properties of $V$.
By Proposition 4.7.4 in \cite{FHL}, any irreducible $V_{\sqrt2E_8}^+\otimes V_{\BW}^+$-module is isomorphic to $M_1\otimes M_2$ for some irreducible $V_{\sqrt2E_8}^+$-module $M_1$ and irreducible $V_{\BW}^+$-module $M_2$.
By Theorem 2.8 in \cite{ADL}, the following fusion rules hold:
\begin{eqnarray}
(M_1\otimes M_4)\times(M_2\otimes M_5)= M_3\otimes M_6\label{Eq:FuTe}
\end{eqnarray}
where $M_1$, $M_2$ and $M_4$, $M_5$ are irreducible modules for $V_{\sqrt2E_8}^+$ and $V_{\BW}^+$ respectively, and $M_1\times M_2=M_3$ and $M_4\times M_5=M_6$.
Since $U^i$ is $2$-elementary totally even, $V$ satisfies the conditions (a)-(c) in Section \ref{Ss} by Lemma \ref{LEx}.
Moreover, $V_{U^k}^+$ is rational since $V_{U^k}^+$ is a framed VOA.
Hence, by Proposition \ref{LScE} $V^\natural$ is a simple current extension of $V$.
More precisely, $V^\natural$ is decomposed into irreducible $V$-modules as follows:

\bl\label{LD3} The $V$-modules $V_\Lambda^+$ and $V_\Lambda^{T,-}$ decompose into irreducible $V$-modules as follows:
\begin{eqnarray*}
V_\Lambda^+\cong\bigoplus_{\lambda+U\in\Lambda/U}(V_{\lambda_1+U^1}^+\otimes V_{\lambda_2+U^2}^+\oplus V_{\lambda_1+U^1}^-\otimes V_{\lambda_2+U^2}^-).\\
V_\Lambda^{T,-}\cong\bigoplus_{\lambda+U\in \Lambda/U}( V_{U^1}^{T_{\chi_{\lambda_1}},+}\otimes V_{U^2}^{T_{\chi_{\lambda_2}},-}\oplus V_{U^1}^{T_{\chi_{\lambda_1}},-}\otimes V_{U^2}^{T_{\chi_{\lambda_2}},+})
\end{eqnarray*}
where $\chi_{\lambda_k}$ is the irreducible character of $\langle(U^k)^*,2\Lambda\rangle/2\Lambda$ defined by $\chi_{\lambda_k}(\gamma)=(-1)^{\langle\lambda_k,\gamma\rangle}$ for $\gamma\in U^k$.
\el
\proof\ Applying Proposition \ref{LD1} and \ref{LD2} to $V^\natural$ and $V_U^+$, we obtain the decomposition of $V^\natural$ into irreducible $V_U^+$-modules.
Let $\theta_{V_L^{T}}$ denote the involution of $V_L^T$ induced by $\theta_{V_L}$.
Since $\theta_{V_\Lambda}=\theta_{V_{U^1}}\otimes\theta_{V_{U^2}}$ and $\theta_{V_\Lambda^{T}}=\theta_{V_{U^1}^T}\otimes\theta_{V_{U^2}^T}$ we obtain the desired decomposition.\qe

Let $S_{10}$ be the set of isomorphism classes of irreducible $V$-modules which appear in Lemma \ref{LD3}.
Then by Proposition \ref{PF} and (\ref{Eq:FuTe}) $S_{10}$ forms an elementary abelian $2$-group of order $2^{10}$ under the fusion rules.
So, we obtain an elementary abelian $2$- subgroup $S_{10}^*$ of $G$.

We define a quadratic form $q$ on $S_{10}$ by setting $q(W)=q_{U^1}(W_1)$ for $W=W_1\otimes W_2\in S_{10}$, where $q_{U^1}$ is a quadratic form given in Proposition \ref{PQF} and $W_k$ is the isomorphism class of an irreducible $V_{U^k}^+$-module.
We note that $q_{U^1}(W_1)=q_{U^2}(W_2)$ for $W=W_1\otimes W_2\in S_{10}$.
Counting singular elements of $S_{10}$, we see that the type of $q$ is plus.
Since the value of $q$ depends only on the graded dimensions, the action of $N_G(S_{10}^*)$ on $S_{10}$ preserves the quadratic form $q$.
Hence $N_G(S_{10}^*)/C_G(S_{10}^*)\subset O^+(10,2)$.

Let $\Phi_{S_{10}}^N$ and $\Phi_{S_{10}}^C$ be the restriction homomorphisms from $N_G(S_{10}^*)$ to $H_{S_{10}}^N$ and from $C_G(S_{10}^*)$ to $H_{S_{10}}^C$, where 
\begin{eqnarray*}
H_{S_{10}}^N&=&\{g\in\Aut(V_{\sqrt2E_8}^+\otimes V_{\BW}^+)|\ S_{10}\circ g= S_{10}\},\\
H_{S_{10}}^C&=&\{g\in\Aut(V_{\sqrt2E_8}^+\otimes V_{\BW}^+)|\ W\circ g= W\ {\rm for\ all}\ W\in S_{10}\}.
\end{eqnarray*}
By Proposition \ref{PSC}, the homomorphisms $\Phi_{S_{10}}^N$ and $\Phi_{S_{10}}^C$ are surjective and $\Ker\Phi_{S_{10}}^N=\Ker\Phi_{S_{10}}^C= S_{10}^*$.
Moreover, $H_{S_{10}}^N/H_{S_{10}}^C\cong N_G(S_{10}^*)/C_G(S_{10}^*)$.
The following lemma will be needed to determine $H_{S_{10}}^N$ and $H_{S_{10}}^C$.

\bl\label{LM3} The automorphism group of $V_{U^1}^+\otimes V_{U^2}^+$ is isomorphic to $\Aut(V_{U^1}^+)\times \Aut(V_{U^2}^+)$.
In particular, its shape is $O^+({10},2)\times2^{16}\cdot \Omega^+(10,2)$.
\el
\proof\ We set $\tilde{H}=\Aut(V_{U^1}^+\otimes V_{U^2}^+)$ and $H=\Aut(V_{U^1}^+)\times \Aut(V_{U^2}^+)$.
Clearly $H\subset \tilde{H}$.

First we consider the stabilizer $K$ of the isomorphism class of $V_{U^1}^-\otimes V_{U^2}^-$ in $\tilde{H}$.
It is easy to see that $\{\lambda+U\in U^*/U|\ \#(\lambda+U)_2=48\}=\phi$.
By Proposition \ref{PSh} (\ref{PSh2}), we have $\Aut(V_U^+)=\AutO(\hat{U})/\langle\theta_{V_U}\rangle$.
Since $V_U^+=V_{U_1}^+\otimes V_{U_2}^+\oplus V_{U_1}^-\otimes V_{U_2}^-$, we have $K=\AutO(\hat{U})/\langle\theta_{V_{U^1}},\theta_{V_{U^2}}\rangle$ by Proposition \ref{PSC}.
By Lemma \ref{LL4}, $K\subset H$.

Next we determine the orbit $Q$ of $V_{U^1}^-\otimes V_{U^2}^-$ under the action of $\tilde{H}$.
Let $S_{{U^k}}$ be the set of all isomorphism classes of irreducible $V_{U^k}^+$-modules.
By Proposition \ref{PSh} (\ref{PSh5}), $S_{U^k}$ is decomposed into $3$-orbits $\{[0]^+\},Q_{U^k},R_{U^k}$ under the action of $\Aut(V_{U^k}^+)$.
We assume that $V_{U^k}^-$ belongs to $Q_{U^k}$.
It is easy to check that
\begin{eqnarray*}
\ch(V_{U^1}^-)\in&8q+q^2\Z[[q]],\ \ch(W_1)&\in q^{1/2}+q^{3/2}\Z[[q]],\\
\ch(V_{U^2}^-)\in&16q+q^2\Z[[q]],\ \ch(W_2)&\in256q^{3/2}+q^{5/2}\Z[[q]],
\end{eqnarray*}
where $W_k\in R_{U^k}$.
Comparing the graded dimensions, we have $Q\subset\{W_1\otimes W_2|\ W_k\in Q_{U^k}\}$.
Since $\Aut(V_{U^k}^+)$ acts transitively on $Q_{U^k}$, we have $Q=\{W_1\otimes W_2|\ W_k\in Q_{U^k}\}$.
In particular, $|Q|=|Q_{U^1}|\times|Q_{U^2}|=527^2$ and $|\tilde{H}:K|=527^2$.

On the other hand, we see that $|H:K|=527^2$.
Thus $\tilde{H}=H$.
\qe

\bt\label{T2} The normalizer $N_G(S_{10}^*)$ is of shape $2^{10+16}\cdot \Omega^+(10,2)$ and the centralizer $C_G(S_{10}^*)$ is of shape $2^{10+16}$.
\et

\proof\ By Proposition \ref{PSh} (\ref{PSh5}) and Lemma \ref{LM3}, $H_{S_{10}}^C\cong 2^{16}$.
Hence $C_G(S_{10}^*)\cong 2^{10+16}$.
Let us determine $H_{S_{10}}^N$.
Let $S_{U^k}$ denote the set of all isomorphism classes of irreducible $V_{U^k}^+$-modules.
Then any element of $S_{U^k}$ ($k=1,2$) appears only once in $S_{10}$.
So, the map $S_{U^1}\to S_{U^2}$, $W_1\mapsto W_2$ ($W_1\otimes W_2\in S_{10}$) is well-defined and bijective.
Let $g\in \Aut(V_{U^1}^+\otimes V_{U^2}^+)$.
By Lemma \ref{LM3}, $g=g_1g_2$ for some $g_k\in \Aut(V_{U^k})$.
Then $g$ is in $H_{S_{10}}^N$ if and only if the actions of $g_k$ on $S_{U^k}$ ($k=1,2$) are the same with respect to the bijection above.
This shows that $H_{S_{10}}^N/H_{S_{10}}^C\cong\Omega^+({10},2)$ by Lemma \ref{LM3}.
Hence $N_G(S_{10})\cong 2^{10+16}\cdot \Omega^+(10,2)$.\qe

\bn\label{NG} {\rm In the previous and this section, we have described five normalizers of elementary abelian $2$-subgroups of $G$ containing $z_{V^\natural}$.
According to \cite{Me,MS}, these normalizers are only maximal $2$-local subgroups of the Monster such that the center contains $2B$ elements.}
\en

\subsection{Monster amalgam}
In this subsection, we show that the set of the groups $\{N_G(S_1^*),N_G(S_2^*),N_G(S_3^*)\}$ forms a Monster amalgam.

We start by recalling the definition of a Monster amalgam from Chapter 5 in \cite{Iv}.
We refer to Chapter 1 of \cite{Iv} for the definition of amalgams.
Let $\mathcal{M}=\{G_1,G_2,G_3\}$ be an amalgam of rank $3$, put $Q_i=O_2(G_i)$, $\bar{G}_i=G_i/Q_i$, $G_{ij}=G_i\cap G_j$, $G_{123}=G_1\cap G_2\cap G_3$, $Z_1=Z(Q_1)$.
Then $\mathcal{M}$ is a {\it Monster amalgam} if the following hold:
\begin{enumerate}[(i)]
\item $Q_1$ is an extraspecial group of order $2^{25}$;
\item $\bar{G}_1\cong Co_1$ acts on $Q_1/Z_1$ as it acts on the Leech lattice $\bar{\Lambda}$ modulo $2$;
\item $G_{123}$ contains a Sylow $2$-subgroup of $G_1$;
\item $\bar{G}_2\cong Sym_{3}\times M_{24}$ and $\bar{G}_3\cong L_3(2)\times 3\cdot Sym_6$;
\item $[G_2:G_{12}]=3$; $[G_3:G_{23}]=[G_{3}:G_{13}]=7$; $[G_3:G_{123}]=21$;
\item for $1\le i<j\le 3$, we have $Q_i\cap Q_j\neq Q_i$.
\end{enumerate}

Set $H_i=N_G(S_i^*)$.
Let $\mathcal{H}$ denote the amalgam $\{H_1,H_2,H_3\}$.

\bt The set $\mathcal{H}$ forms a Monster amalgam.
In particular the group $\langle H_i|\ i=1,2,3\rangle$ is isomorphic to the Monster simple group.
\et
\proof\ 
By using the description of the shapes of $H_i$ in Theorem \ref{T1}, we will show that $\mathcal{H}$ satisfies the conditions (i)-(vi).
Set $R_i=O_2(H_i)$, $\bar{H}_i=H_i/R_i$, $H_{ij}=H_i\cap H_j$, $H_{123}=H_1\cap H_2\cap H_3$ and $Z=Z(H_1)$.

First we consider the group $H_1$.
Since $H_1/S_1^*$ is the automorphism group of $V_\Lambda^+$, $\bar{H}_1\cong Co_1$ acts on $H_1/S^*_1$ as it acts on the Leech lattice $\bar{\Lambda}$ modulo $2$.
In particular, this action is irreducible.
Since the shape of $R_1$ is $2.2^{24}$, $\Cent(R_1)=S^*_1$ or $R_1$.
On the other hand, in \cite{Sh1} an automorphism in $R_1$ of order $4$ commuting with $S^*_1$ is described in terms of $V_L^+$.
This shows that $R_1$ is an extraspecial $2$-group, and (i) follows.
Since $Z(R_1)=S^*_1$, (ii) also follows.

Next, we consider $H_{ij}$ and $H_{123}$.
The inclusion $S_1^*\subset S_2^*$ shows that $C_{G}(S_2^*)\subset H_{12}$.
Recall that $H_{2}/C_{G}(S_2^*)=GL(S_2^*)\cong L_2(2)\cong Sym_3$.
Since $H_{12}$ is the stabilizer of $S^*_1$ in the action of $H_2$ on $S_2^*$, we obtain $H_{12}/C_{G}(S_2^*)\cong\Z_2$.
Hence we obtain $[H_2:H_{12}]=3$.
Similarly, for $i=1,2$ the inclusion $S_i^*\subset S_3^*$ shows that $C_{G}(S_3^*)\subset H_{i3}$.
Recall that $H_{3}/C_{G}(S_3^*)=GL(S_3^*)\cong L_3(2)$.
For $i=1,2$ the group $H_{i3}$ is the stabilizer of $S_i^*$ in $H_3$, which shows that $[H_3:H_{i3}]=7$.
Since $H_{123}$ is the stabilizer of both $S_1^*$ and $S_3^*$ in $H_3$, we obtain $[H_3:H_{123}]=21$.
Hence (v) follows.
Moreover the shape of $H_{123}$ is $(2^3.2^{20}.(2^{4+12}.3Sym_6)).2^3$.
In particular a Sylow subgroup of $H_{123}$ is also one of $H_1$, and (iii) holds.

By the shape of $H_i$, we obtain $|R_1|=2^{25}$, $|R_2|=2^{35}$ and $|R_3|=2^{39}$.
Moreover $R_1\cap R_2\cong 2^{2}.2^{22}$, $R_1\cap R_3\cong 2^{3}.2^{20}$ and $R_2\cap R_3\cong 2^3.2^{20}.2^{10}$.
Hence (vi) follows.

Finally let us consider $\bar{H}_i$.
By the shape of $H_i$, we obtain $\bar{H}_2\cong M_{24}.Sym_3$ and $\bar{H}_3\cong 3Sym_6.L_3(2)$.
Let $K_2$ denote the subgroup of all elements of $\bar{H}_2$ acting trivially on $M_{24}$.
Since the Schur multiplicity of $M_{24}$ is $1$, $K_2$ is isomorphic to $Sym_3$ and is compliment to the subgroup of shape $M_{24}$.
Hence $\bar{H}_2\cong M_{24}\times Sym_3$.
Let $K_3$ denote the subgroup of all elements of $\bar{H}_3$ acting trivially on $3Sym_6$ and let $I$ denote the normal subgroup of $\bar{H}_3$ of shape $3Sym_6$.
Then $K_3\cap I=Z(I)$.
Since $K_3$ is normal in $\bar{H}_3$, so is $IK_3/I$ in $\bar{H}_3/I$.
By Lemma \ref{LL2} and (\ref{Eq:SeqNC}), $K_3$ contains a subgroup of shape $Sym_3$.
Since $L_3(2)$ is simple, $Z(I)K_3/Z(I)$ is isomorphic to $L_3(2)$ and is compliment to the subgroup of shape $3Sym_6$.
Hence $\bar{H}_3\cong L_3(2)\times 3Sym_6$.
Thus (iv) follows.
Therefore $\mathcal{H}$ is a Monster amalgam.
By Proposition 5.15.1 in \cite{Iv} $\langle H_i|\ i=1,2,3\rangle$ is isomorphic to the Monster.\qe

\end{document}